\documentclass[12pt,a4paper]{article}

\usepackage[latin1]{inputenc}
\usepackage{color}

\usepackage{bussproofs}
\usepackage{graphicx}
\usepackage{graphics}
\usepackage{latexsym}
\usepackage{amssymb,amsmath}
\usepackage{euscript}
\usepackage{enumerate}
\usepackage{amsmath}
\usepackage{amsfonts}
\usepackage{multicol}
\usepackage{etex}
\usepackage[all]{xy}
\usepackage{color}

\newtheorem{teo}{Theorem}[section]

\newtheorem{lema}[teo]{Lemma}
\newtheorem{prop}[teo]{Proposition}

\newtheorem{rem}[teo]{Remark}
\newenvironment{dem}{\noindent {\em Proof:} }{\hfill $\square$ \bigskip}
\newenvironment{prueba}[1]{\noindent\bf Proof of #1. \rm}{$\quad \hfill \square$ \bigskip}
\numberwithin{equation}{section}
\newcommand{\R}{\mathbb{R}}

\newcommand{\Z}{\mathbb{Z}}
\newcommand{\N}{\mathbb{N}}

\def\XXint#1#2#3{{\setbox0=\hbox{$#1{#2#3}{\int}$}
    \vcenter{\hbox{$#2#3$}}\kern-.5\wd0}}

\usepackage{fancyhdr}

\fancypagestyle{encabezadotesis}{
\fancyhead{}
\fancyhead[RO]{\em \nouppercase{\rightmark}}
\fancyhead[LE]{\bf \nouppercase{\leftmark}}
\cfoot{\thepage}
}

\begin{document}
\title{Two weighted estimates for  generalized fractional maximal operators on non homogeneous spaces.}
\author{Gladis Pradolini and Jorgelina Recchi}
\date{}
\maketitle

\renewcommand{\thefootnote}{\fnsymbol{footnote}}
\footnotetext{2010 {\em Mathematics Subject Classification}:
42B25} \footnotetext {{\em Key words and phrases}:
non-homogeneous space, generalized fractional operators, weights} \footnotetext{
The first author is supported by Instituto de Matem\'atica Aplicada del Litoral (CONICET-UNL) and Departamento de Matem\'atica
(FIQ-UNL), Santa Fe, Argentina.
}
\footnotetext{The second author is  supported by Instituto de Matem\'atica Bah\'ia Blanca (CONICET-UNS) and Departamento de Matem\'aticas (UNS), Bah\'ia Blanca, Argentina.}
\begin{abstract}
Let $\mu$ be a non-negative Borel measure on $\R^d$ satisfying that
the measure of a cube $\R^d$ is smaller than the length of its side
raised to the $n$-th power, $0<n\leq d$.

In this article we study the class of weights related to the
boundedness of radial fractional type maximal operator associated to
a Young function $B$ in the context of non-homogeneous spaces
related with the measure $\mu$. This type of maximal operators are
the adequate operators related with commutators of singular and
fractional operators. Particularly, we give an improvement of a two
weighted result for certain fractional maximal operator proved in
\cite{wang tan lou}.
\end{abstract}

\section{Introduction and statements of the main results}

Let $\mu$ be a non-negative upper Ahlfors $n$-dimensional measure on
$\R^d$, that is, a Borel measure satisfying
\begin{equation}\label{mu}
\mu(Q)\leq l(Q)^n
\end{equation}
for any cube $Q\subset \R^n$ with sides parallel to the coordinate
axes, where $l(Q)$ stands for the side length of $Q$ and $n$ is a
fixed real number such that $0<n\leq d$. Besides, for $r>0$, $rQ$
will mean the cube with the same centre as $Q$ and with
$l(rQ)=rl(Q)$.

In the last decades, this measure have proved to be adequate for the
development of many results in Harmonic Analysis which were known
that hold in the context of doubling measures, that is, Borel
measures $\nu$ for which there exists a positive constant $D$ such
that $\nu(2Q)\le D \nu(Q)$ for every cube $Q\subset \R^d$. For
example, many interesting results related with different operators
and spaces of functions with non doubling measures can be found in
\cite{nazarov treil volverg 1}, \cite{nazarov treil volverg 2},
\cite{mateu mattila nicolau orobitg}, \cite{tolsa}, \cite{garcia
cuerva martell} and \cite{meng yang} between a vast bibliography on
this topic.

In \cite{wang tan lou} the authors studied two weighted norm
inequalities for a fractional maximal operator associated to a
measure $\mu$ satisfying condition (\ref{mu}). Concretely, they
considered the following version of the fractional maximal operator
defined, for $0\le\alpha<1$, by
\begin{equation*}
{M}_{\alpha}f(x)= \sup_{Q\ni x}\frac{1}{\mu(5Q)^{1-\alpha}}\int_Q
|f(y)|d\mu(y),
\end{equation*}
and proved the following result.
\begin{teo}\label{chinos}Let $1<p<q<\infty$ and $0\le\alpha<1$. Let $(u,v)$ be
a pair of weights such that for every cube $Q$
\begin{equation}\label{muwtl}
l(Q)^{n(1-1/p)}\mu(Q)^{\alpha-1}u(3Q)^{\frac{1}{q}}\|v^{-\frac{1}{p}}\|_{\Phi,Q}\leq
C
\end{equation}
where $\Phi$ is a Young function whose complementary function
$\bar{\Phi}\in B_p$. Then

\begin{equation*}
\left(\int_{\R^d}M_{\alpha}f(x)^q u(x)\, d\mu(x)\right)^{1/q}\le
C\left(\int_{\R^d}|f(x)|^p v(x)\, d\mu(x)\right)^{1/p}.
\end{equation*}
for every $f\in L^p(v)$ bounded with compact support.
\end{teo}
The radial Luxemburg type average in theorem above is defined by
$$\|\cdot\|_{\Phi,Q}=\inf\{\lambda>0:
\frac{1}{l(Q)^n}\int_{Q}\Phi\left(\frac{|f(x)|}{\lambda}\right)\,d\mu(x)\le1\},$$
and a Young function $B$ satisfies the $B_p$ condition,
$1<p<\infty$, if there is a positive constant $c$ such that
\begin{equation*}\label{Bp}
\int_c^{\infty} \frac{B(t)}{t^p}\frac{dt}{t}<\infty.
\end{equation*}

Let us make some comments about Theorem $\ref{chinos}$. When $\mu$
is the Lebesgue measure and $u=v=1$, it is easy to note that
condition ($\ref{muwtl}$) holds if and only if $1/q=1/p-\alpha$ for
any $\Phi$ as in the hypothesis. On the other hand, if we consider
an upper Ahlfors $n$-dimensional measure $\mu$ and if we take
$\Phi(t)=t^{rp'}$, for $1<r<\infty$, $1/q=1/p-\alpha$ and $u=v=1$ in
condition ($\ref{muwtl}$) we have that if the following inequality
holds
$$l(Q)^{n(1-1/p)}\mu(Q)^{\alpha-1}\mu(3Q)^{\frac{1}{p}-\alpha} \left(\frac{\mu(Q)}{l(Q)^n}\right)^{1/rp'}\leq
C,$$ then
$$\left(\frac{l(Q)^n}{\mu(Q)}\right)^{1/(p'r')}\le C$$
which implies that the measure $\mu$ satisfying the growth condition
($\ref{mu}$) also satisfies the ``lower" case, that is $\mu(Q)\ge C
l(Q)^n$, with a constant independent of $Q$. So, the weights $u=v=1$
are not allowed in this case unless the measure is Ahlfors, that is
$\mu(Q)\simeq l(Q)^n$, for every cube $Q$.
%---------------------------------------------------
Moreover, let $$Mu(x)=\sup_{Q\ni x}\frac{1}{\mu(Q)}\int_Q
|u(y)|d\mu(y).$$ When $\mu$ is the Lebesgue measure and
$\Phi(t)=t^{rp'}$, it is easy to check that the pair of weights $(u,
(Mu)^{p/q})$, with $1/q=1/p-\alpha$, satisfies condition
($\ref{muwtl}$). On the other hand, suppose that this pair satisfies
the same condition for a measure satisfying ($\ref{mu}$) and let
$u\in A_1(\mu)$. Thus, the following chains of inequalities holds
\begin{eqnarray*}
C&\ge&
\frac{l(Q)^{n/p'}}{\mu(Q)^{1-\alpha-1/q}}\left(\frac{1}{\mu(Q)}\int_Q
u\, d\mu\right)^{1/q}\left(\frac{1}{l(Q)^n}\int_Q
\left((Mu)^{p/q}\right)^{-rp'/p}\right)^{1/(rp')}\\
&\ge&
\frac{l(Q)^{n/p'-n/(rp')}}{\mu(Q)^{1/p'-1/(rp')}}\left(\frac{1}{\mu(Q)}\int_Q
u^{p/q}\, d\mu\right)^{1/p}\left(\frac{1}{\mu(Q)}\int_Q
\left(u^{p/q}\right)^{-rp'/p}\right)^{1/(rp')}\\
&\ge&
\frac{l(Q)^{n/(r'p')}}{\mu(Q)^{1/(r'p')}}\left(\frac{1}{\mu(Q)}\int_Q
u^{p/q}\, d\mu\right)^{1/p}\left(\frac{1}{\mu(Q)}\int_Q
\left(u^{p/q}\right)^{-p'/p}\right)^{1/(p')}\\
&\ge& \frac{l(Q)^{n/(r'p')}}{\mu(Q)^{1/(r'p')}},
\end{eqnarray*}
which implies again that $\mu$ must be an Ahlfors measure.

In \cite{garcia cuerva martell} the authors considered the radial
fractional maximal function associated to an upper Ahlfors
$n$-dimensional measure $\mu$ which is defined, for $0\leq
\alpha<n$, by
$$
\mathcal{M}_{\alpha}f(x)= \sup_{Q\ni
x}\frac{1}{l(Q)^{n-\alpha}}\int_Q |f(y)|d\mu(y).
$$
In the same article they study weighted boundedness properties for
$\mathcal{M}_{\alpha}$ on non homogeneous spaces.

%------------------------------------------------------

In this paper we introduce a generalized version of the radial
fractional maximal operator defined in \cite{garcia cuerva martell},
associated to a Young function $B$. This type of maximal operators
are not only a generalization but also they have proved to be the
adequate operators related with commutators of singular and
fractional integral operators in different settings, (see for
example \cite{bernardis hartzstein pradolini}, \cite{perez1},
\cite{perez2}, \cite{perez3}, \cite{perez pradolini}, \cite{lorente
martell riveros de la torre}, \cite{lorente riveros de la torre} and
\cite{bernardis lorente riveros}). It is important to note that the
examples of weights given above satisfy the condition obtained in
our theorem when $\mu$ is an upper Ahlfors $n$-dimensional measure.
In this sense, when $B(t)=t$ , our result is better than the
corresponding result in \cite{wang tan lou}.

\medskip

%------------------------------------------------------
In order to state the main results we introduce some preliminaries.
Given a Young function $B$, we define $L^B_\mu(\R^d)$ as the set of
all measurable functions $f$ for which there exists a positive
number $\lambda$ such that
$$
\int_{\R^d}B\left(\frac{|f(x)|}{\lambda}\right)\,d\mu(x)<\infty.
$$

The radial fractional type maximal operator associated to a Young
function $B$ is defined by
$$
\mathcal{M}_{\alpha,B}(f)(x)= \sup_{Q\ni x} \; l(Q)^{\alpha}
\|f\|_{B,Q}, \quad 0\le \alpha<n,
$$
where $$\|f\|_{B,Q}=\inf\{\lambda>0:
\frac{1}{l(Q)^n}\int_{Q}B\left(\frac{|f(x)|}{\lambda}\right)\,d\mu(x)\le1\}$$
is the radial Luxemburg average (see \S\ref{definicion de norma de
lux}). When $B(t)=t$ then
\begin{equation}\label{promedio radial}
\|f\|_{B,Q}=\frac{1}{l(Q)^n}\int_Q |f|\, d\mu.
\end{equation}
When $\alpha=0$, we write $\mathcal{M}_{0,B}=\mathcal{M}_{B}$.

\medskip

The following theorem gives sufficient conditions for strong type
inequalities for $\mathcal{M}_{\alpha,B}$ on non homogeneous spaces.

\begin{teo}\label{teo similar al 5.3 de GC martell}
Let $1<p<q<\infty$, $0\leq\alpha<n$ and let $\mu$ be an upper
Ahlfors $n$-dimensional measure in $\R^d$. Let $B$ be a
submultiplicative Young function  such that $B^{q_0/p_0}\in B_{q_0}$
for some $1<p_0\le n/\alpha$ and $1/q_0=1/p_0-\alpha/n$, and let
$\phi$ and $\varphi$ be two Young functions such that $C_1
\varphi^{-1}(t)t^{\alpha/n}\leq \;B^{-1}(t)\le C_2
\phi^{-1}(t)t^{\alpha/n}$ for some positive constants $C_1$ and
$C_2$. If $A$ and $C$ are two Young functions such that
$A^{-1}C^{-1}\preceq B^{-1}$ with $C \in B_p$ and $(u,v)$ is a pair
of weights such that for every cube $Q$
\begin{equation}\label{good}
l(Q)^{\alpha-\frac{n}{p}}u(3Q)^{\frac{1}{q}}\|v^{-\frac{1}{p}}\|_{A,Q}\leq
K
\end{equation}
then, for all $f\in L_\mu^p(v)$.
$$
\|\mathcal{M}_{\alpha,B}(f)\|_{L^q_\mu(u)}\leq C\; \|f\|_{L^p_\mu(v)},
$$
\end{teo}

\medskip
\begin{rem}When $u=v=1$ and $1/q=1/p-\alpha/n$ then condition
$(\ref{good})$ is satisfied for any upper Ahlfors $n$-dimensional
measure $\mu$. Thus, this result is an improvement of that given in
\cite{wang tan lou} in the sense that the unweighted boundedness of
the operator is obtained for any measure satisfying the growh
condition $(\ref{mu})$. The same is true for the second example
considered above.
\end{rem}

\begin{rem}\rm \label{remark B}
When $B(t)=t \log(e+t)^k$ it can be easily seen that $B$ is
submultiplicative, $B^{q_0/p_0}\in B_{q_0}$ for every $p_0, q_0
>1$ and $$
B^{-1}(t)\approx t^{\alpha/n} \frac{t^{1-\alpha/n}}{\log(e+t)^k}\approx t^{\alpha/n}\phi^{-1}(t),
$$
when $\phi(t)=\left(t\log(e+t)^k\right)^{\frac{n}{n-\alpha}}$.
Moreover, the functions $A(t)=t^{rp'}$ and
$C(t)=(t\log(e+t)^k)^{(rp')'}$ satisfy $$A^{-1}C^{-1}\preceq
B^{-1}.$$ For $\delta>0$, other examples are given by
$A(t)=t^{p'}\log(e+t)^{(k+1)p'-1+\delta}$ and
$C(t)=t^{p}\log(e+t)^{-(1+\delta(p-1))}$, (see \cite{cruz uribe
perez}).

\bigskip

It is also important to note that Theorem 3.1 in \cite{garcia cuerva
martell} is a special case of the previous theorem by considering
$A(t)=t^{rp'}$, $C(t)=t^{(rp')'}$ and $B(t)=t$.
\end{rem}

Let us make some comments about the upper Ahlfors $n$-dimensional
measure $\mu$ satisfying (\ref{mu}). It is well known that for such
measures the Lebesgue differentiation theorem holds; that is, for
every $f\in L_{loc}^1(\R^d)$ and a.e. $x$
$$
\frac{1}{\mu(Q)}\int_Q f(y)\,d\mu(y)\rightarrow f(x),
$$
when $Q$ decreases to $x$ (see \cite{yang yang hu}). However, if we
take radial averages like those defined in (\ref{promedio radial})
this is not longer true. In fact, let us consider $\mu$ defined  by
$d\mu(t)=e^{-t^2}dt$, which is an upper Ahlfors $1$-dimensional
measure, and $f(t)=e^{\theta t^2}$, $\theta \in \R$. Let $x \in \R$,
then
$$
\lim_{r\to 0}\frac{1}{2r}\int_{x-r}^{x+r}f(t)d\mu(t)=\lim_{r\to
0}\frac{1}{2r}\int_{x-r}^{x+r}e^{(\theta-1)t^2}\,dt=e^{(\theta-1)x^2}
$$
which differs from $f$ in a.e. $x$.

Given a Young function $B$, let $h_B$ be the function defined by
\begin{equation*}
h_B(s)= \sup_{t>0}\frac{B(st)}{B(t)}, \;\; 0\leq s<\infty.
\end{equation*}
If $B$ is submultiplicative then $h_B \approx B$. More generally,
given any $B$, for every $s$, $t\geq0$, $B(st)\leq h_B(s)B(t)$. It
is easy to proof (see \cite[Lemma 3.11]{cruz uribe fiorenza 2002}),
that if $B$ is a Young function then $h_B$ is nonnegative,
submultiplicative, increasing in $[0,\infty)$, strictly increasing
in $[0,1]$ and $h_B(1)=1$.

\medskip

The following theorem gives an modular endpoint estimate for
$\mathcal{M}_{\alpha,B}$ on non-homogeneous spaces. When $\mu$ is
the Lebesgue measure, this result was proved in \cite{cruz uribe
fiorenza}.

\begin{teo}\label{tipo debil de la malphaB}
Let $0\le\alpha<n$ and let $\mu$ an upper Ahlfors $n$-dimensional
measure on $\R^d$. Let $B$ be a Young function and suppose that, if
$\alpha>0$,  $B(t)/t^{\frac{n}{\alpha}}$ is decreasing for all
$t>0$. Then there exists a constant $C$ depending only on $B$ such
that for all $t>0$, $\mathcal{M}_{\alpha,B}$ satisfies the following
modular weak-type inequality
$$
\phi\left[\mu \left(\{x \in \R^d: \mathcal{M}_{\alpha,B}(f)(x)>t\}\right)\right]\leq \;C\;\int_{\R^n} B\left(\frac{|f(y)|}{t}\right)d\mu(y),
$$
for all $f\in L_\mu^B(\R^d)$, where $\phi$ is any function such that
$$
\phi(s)\leq C_1 \phi_1(s)= \left\{
                     \begin{array}{ll}
                       0 & \text {if }  s=0  \\ \\
                       \frac{s}{h_B(s^{\alpha/n})} & \text{if } s>0\\ \\
                     \end{array}
                   \right.
$$
\end{teo}

\begin{rem}\rm
It is easy to see that the function $B(t)=t\log(e+t)$ satisfies the
hypothesis of the theorem above and thus
 $$
\mu \left(\{x \in \R^n: \mathcal{M}_{\alpha,B}(f)(x)>t\}\right)\leq \;C\;\psi\left[\int_{\R^n} B\left(\frac{|f(y)|}{t}\right)d\mu(y)\right],
$$
for all $f\in L_\mu^B(\R^d)$, where $\psi=[t
\log(e+t^{\alpha/n})]^{n/(n-\alpha)}$. In the context of spaces of
homogeneous type this last result was proved in \cite{gorosito
pradolini salinas}.
\end{rem}
The following result is a pointwise estimate between the radial
fractional type maximal operator associated to a Young function $B$
and the radial maximal operator associated to a Young $\psi$ related
with $B$ on non homogeneous spaces.

\begin{teo}\label{teo puntual de Malpha a M}
Let $0\leq \alpha < n$ and $1<p<n/\alpha$. Let $\mu$ be an upper
Ahlfors $n$-dimensional measure. Let $q$ and $s$ be defined by
$1/q=1/p-\alpha/n$ and $s=1+q/{p'}$, respectively. Let $B$ and
$\phi$ be Young functions such that $\phi^{-1}(t)t^{\alpha/n}\geq
C\; B^{-1}(t)$ and $\psi(t)=\phi(t^{1-\alpha/n})$. Then for every
measurable function $f$, the following inequality
$$
\mathcal{M}_{\alpha,B}(f)(x)\leq C\;
\mathcal{M}_{\psi}(|f|^{p/s})(x)^{1-\alpha/n}\left(\int_{\R^d}
|f(y)|^pd\mu(y)\right)^{\alpha/n}
$$
holds in a.e. $x$.
\end{teo}

When $\mu$ is the Lebesgue measure, the result above was proved in
\cite{bernardis dalmasso pradolini} (see also \cite{pradolini2} for
similar multilinear versions and \cite{gorosito pradolini salinas 2}
for the case $B(t)=t$, both in the euclidean context).

\bigskip
The next theorem gives sufficient conditions on the function $B$ in
order to obtain the boundedness of $\mathcal{M}_B$ on $L^p(\mu)$. In
the euclidean context, this result was proved in \cite{perez On
suff} and in \cite{pradolini salinas} in the framework of spaces of
homogeneous type.

\begin{teo}\label{Meta}
Let $\mu$ be an upper Ahlfors $n$-dimensional measure. Let $B$ be a
Young function such that $B \in B_p$, then,
$$
\mathcal{M}_B :L^p(\mu)\to L^p(\mu).
$$
\end{teo}

The following result is a fractional version of Theorem \ref{Meta}
and gives a sufficient condition on the function $B$ that guarantees
the continuity of the radial fractional type maximal operator
$\mathcal{M}_{\alpha,B}$ between Lebesgue spaces with non necessary
doubling measures.

\begin{teo}\label{acotmaximal}
Let $\mu$ be an upper Ahlfors $n$-dimensional measure. Let $0<\alpha
< n$ and $1<p\leq n/\alpha$. Let $q$ and $s$ be defined by
$1/q=1/p-\alpha/n$ and $s=1+q/{p'}$ respectively. Let $B$ be a
submultiplicative Young function such that $B^{q/p}\in B_q$ and let
$\phi$ be a Young function such that $\phi^{-1}(t)t^{\alpha/n}\geq
C\;B^{-1}(t)$. Then
$$
\mathcal{M}_{\alpha,B}:L^p(\mu)\to L^q(\mu).
$$
\end{teo}

The next theorem is  very interesting since it allows us to readily find examples of $A_1$ weights on  non homogeneous spaces.

 \begin{teo} \label{malpha es un peso A1}
 Given $\alpha$, $0<\alpha<n$, and a non-negative function $f$. There exists a constant $C$ such that,
 $$
 \mathcal{M}(\mathcal{M}_\alpha f)(x)\leq \; C\;\mathcal{M}_\alpha f(x).
 $$
 \end{teo}
 \begin{dem}
 Fix a cube $Q$. We shall see that,
 $$
 \frac{1}{l(Q)^n}\int_Q \mathcal{M}_\alpha f(y)\;d\mu(y)\leq \;C\;\mathcal{M}_\alpha f(x)\;\;\;\;\; \text{for a.e.}\;\;\;\; x\in Q.
 $$
 with $C$ independent of $Q$. Let $\widetilde{Q}=Q^3$, the 3-dilated of  $Q$. We write $f=f_1+f_2$ with $f_1=f\chi_{\widetilde{Q}}$. Then, $\mathcal{M}_\alpha f(x)\leq \mathcal{M}_\alpha f_1(x)+\mathcal{M}_\alpha f_2 (x)$.
 $$
 \frac{1}{l(Q)^n}\int_Q \mathcal{M}_\alpha f_1 (y)d\mu(y)=\frac{1}{l(Q)^n}\int_0^\infty \mu \{x\in Q: \mathcal{M}_\alpha f_1(x)>t\}dt\leq
 $$
 $$
 \leq \frac{1}{l(Q)^n}\left(\mu(Q)R+\int_R^\infty \mu \{x\in Q: \mathcal{M}_\alpha f_1(x)>t\}dt\right)
 $$
 By \cite[Prop 2.1]{garcia cuerva martell}, we know that $\|\mathcal{M}_\alpha f\|_{L^{\frac{n}{n-\alpha},\infty}(\mu)}\leq \|f\|_{L^{1}(\mu)}$. Then, since $\mu(Q)\leq l(Q)^n$
 $$
 \frac{1}{l(Q)^n}\int_Q \mathcal{M}_\alpha f_1 (y)d\mu(y)\leq R+ \frac{c}{l(Q)^n} \|f_1\|_{L^1(\mu)}^{n/{n-\alpha}}\int_R^\infty \frac{dt}{t^{n/{n-\alpha}}}.
 $$
  By taking $R=\frac{\|f_1\|_{L^1(\mu)}}{l(Q)^{n-\alpha}}$, we get
 $$
  \frac{1}{l(Q)^n}\int_Q \mathcal{M}_\alpha f_1 (y)d\mu(y)\leq \;C_{\alpha,n}\; \frac{\|f_1\|_{L^1(\mu)}}{l(Q)^{n-\alpha}}=  \frac{C_{\alpha,n}}{l(\widetilde{Q})^{n-\alpha}}\int_{\widetilde{Q}}f(y)d\mu(y)\leq C_{\alpha,n}\mathcal{M}_\alpha f(x)
 $$
 for every $x \in Q$.

 To deal with $\mathcal{M}_\alpha f_2$ is even simpler. It is enough to see that, because of the fact that $f_2$ lives far from $Q$ (outside $\widetilde{Q}$), for any two points $x,y \in Q$, we have $\mathcal{M}_\alpha f_2 (x)\leq C \mathcal{M}_\alpha f_2(y)$, with $C$ an absolute constant. Indeed if $Q_0$ is a cube containing $x$ and meeting  $\R^{n}\setminus \widetilde{Q}$, then $Q\subset Q_0^3$, so that
 $$
 \frac{1}{l(Q_0)^{n-\alpha}}\int_{Q_0} f_2(t)d\mu(t)\leq \frac{3^{n-\alpha}}{l(Q_0^3)^{n-\alpha}}\int_{Q_0^3}f_2(t)d\mu(t) \leq 3^{n-\alpha}\mathcal{M}_\alpha f_2(y).
 $$
 Thus
 $$
 \frac{1}{l(Q)^n}\int_Q \mathcal{M}_\alpha f_2 (y)d\mu(y)\leq\; C \frac{\mu(Q)}{l(Q)^n} \mathcal{M}_\alpha f(x)\leq\; C \mathcal{M}_\alpha f(x)
 $$
 for every $x\in Q$.
 \end{dem}

 \section{Preliminaries and auxiliar lemmas}

A function $B:[0,\infty)\to [0,\infty)$ is a Young function if it is
convex and increasing, if $B(0)=0$, and if $B(t) \to \infty$ as $t
\to \infty$. We also deal with submultiplicative Young functions,
which means that $B(st)\le B(s)B(t)$ for every $s$, $t>0$. If $B$ is
a submultiplicative Young function, it follows that $B'(t)\simeq
B(t)/t$ for every $t>0$

Given a Young function $B$ and a cube $Q$, we define the radial
Luxemburg average of $f$ on $Q$ associated to $\mu$ by
\begin{equation} \label{definicion de norma de lux}
\|f\|_{B,Q}=\inf\left\{\lambda>0: \frac{1}{l(Q)^n}\int_Q
B\left(\frac{|f(x)|}{\lambda}\right)d\mu(x)\le 1\right\}.
\end{equation}

The radial Luxemburg average has two rescaling properties which we
will use repeatedly. Given any Young function $A$ and $r>0$,
$$
\|f^r\|_{A,Q}=\|f\|_{B,Q}^{r},
$$
where $B(t)=A(t^r)$. By convexity, if $\tau>1$

$$
\|f\|_{A,Q}\leq \tau^n \|f\|_{A,\tau Q}.
$$

Given a Young function $B$, the complementary Young function
$\tilde{B}$ is defined by
$$
\tilde{B}(t)=\sup_{s>0}\{st-B(s)\},\;\;\; t>0.
$$It is well known that $B$ and $\tilde{B}$ satisfy the following
inequality
$$t\le B^{-1}(t)\tilde{B}^{-1}\le 2t.$$

It is also easy to check that the following version on the
H\"{o}lder's inequality
\begin{equation*}
\frac{1}{l(Q)^n}\int_{Q} |f(x)g(x)|d\mu(x)\leq
2\|f\|_{B,Q}\|g\|_{\tilde{B},Q}
\end{equation*}
holds. Moreover, there is a further generalization of the inequality
above. If $A$, $B$ and $C$ are Young functions such that for every
$t\geq t_0>0$,
$$
B^{-1}(t)C^{-1}(t)\leq c\, A^{-1}(t),
$$
then, the inequality
\begin{equation}\label{desigualdad de holder general}
\|fg\|_{A,Q}\leq K \|f\|_{B,Q}\|g\|_{C,Q}
\end{equation}
holds.

In this papers we give boundedness results for the maximal operator
$\mathcal{M}_{\alpha,B}$ between Lebesgue spaces. We begin with an
usefull property related with $B_p$ condition.

\begin{prop}\label{Bq}
Let $B$ be a submultiplicative Young function and let $\phi$ be a
Young function such that $\phi^{-1}(t)t^{\alpha/n}\geq
C\,B^{-1}(t)$. Let $1<p<n/\alpha$, $1/q=1/p-\alpha/n$ and
$s=q(1-\alpha/n)$. If $B^{q/p}\in B_q$, then the function $\psi$
defined by $\psi(t)=\phi(t^{1-\alpha/n})$ belongs to $B_s$.
\end{prop}

\begin{dem} From the definition of $\psi$ and by changing variables
we obtain that
$$
\int_1^\infty \frac{\psi(t)}{t^s}\frac{dt}{t}= \int_1^{\infty}
\frac{\phi(t^{1-\alpha/n})}{t^s}\frac{dt}{t}
=\left(\frac{n}{n-\alpha}\right)\int_1^{\infty}
\frac{\phi(r)}{r^{ns/{(n-\alpha)}}}\frac{dr}{r}.
$$
From the relation between $B$ and $\phi$ it is easy to see that
$\phi$ is a submultiplicative function. Thus, noting that
$q=ns/{(n-\alpha)}$ we obtain
\begin{eqnarray*}
\int_1^{\infty} \frac{\phi(r)}{r^{ns/{(n-\alpha)}}}\frac{dr}{r}&=&
\int_1^{\infty} \frac{\phi(r)}{r^{q}}\frac{dr}{r}\\&\leq& c\,
\int_c^{\infty}\frac{u^{1+q\alpha/n}}{\left(\phi^{-1}(u)u^{\alpha/n}\right)^q}\frac{du}{u}\\
&\leq& C\, \int_c^{\infty} \frac{u^{q/p}}{B^{-1}(u)^q}\frac{du}{u}=
C\, \int_c^{\infty}\frac{B(t)^{q/p}}{t^q}\frac{dt}{t}<\infty.
\end{eqnarray*}
\end{dem}

The proof of Theorem \ref{tipo debil de la malphaB} requires any
lemmas. The first of them was proved in \cite{cruz uribe fiorenza}
and the second in \cite{hardy littlewood polya}.

\begin{lema}\label{el 3.12}
Given $0\leq\alpha<n$, let $B$ be a Young function such that for
$\alpha > 0$, $B(t)/t^{n/\alpha}$ is decreasing for all $t>0$. Then
the function $\phi_1$ en Theorem \ref{tipo debil de la malphaB} is
increasing, and $\phi_1(s)/s$ is decreasing. Moreover, there exists
$\phi$ such that $\phi(s)\leq C_1 \phi_1(s)$ and $\phi$ is
invertible.
\end{lema}

\begin{lema}\label{el 3.13}
If $\phi(t)/t$ is decreasing, then for any positive sequence $\{x_j\}$,
$$
\phi\left(\sum_j x_j\right)\leq \sum_j\phi(x_j).
$$
\end{lema}

The following lemma is a generalization of Lemma 3.2 in \cite{garcia
cuerva martell} for radial Luxemburg type averages. When $\mu$ is
the Lebesgue measure, it was proved in \cite{cruz uribe fiorenza}.

\begin{lema} \label{paso de un cubo a un diadico}
Suppose that $0\le \alpha <n$, $B$ is a Young function and $f$ is a
nonnegative bounded function with compact support. If for $t>0$ and
any cube $Q$
$$
l(Q)^{\alpha}\|f\|_{B,Q}>t,
$$
then, there exist a dyadic cube $P$ such that $Q\subset 3P$ satisfying
$$
l(P)^{\alpha}\|f\|_{B,P}>\beta t,
$$
where $\beta$ is a nonnegative constant.
\end{lema}

\begin{dem}Let $Q$ be a cube with
\begin{equation}\label{para un cubo cualquiera}
l(Q)^{\alpha}\|f\|_{B,Q}>t.
\end{equation}
Let $k$ be the unique integer such that $2^{-(k+1)}<l(Q)\leq
2^{-k}$. There is some dyadic cube with side length $2^{-k}$, and at
most $2^d$ of them, $\{J_i: i=1,...,N\}$, $N\le 2^d$, meeting the
interior of $Q$. It is easy to see that for one of these cubes, say
$J_1$,
$$
\frac{t}{2^d}<l(Q)^\alpha\|\chi_{J_1}f\|_{B,Q}.
$$
This can be seen as follows. If for each $i=1,2,.., N$ we have
$$
l(Q)^{\alpha}\|\chi_{J_i}f\|_{B,Q}\leq \frac{t}{2^d},
$$
since $Q\subset \cup_{i=1}^{N}J_i$ we obtain that
\begin{eqnarray*}
l(Q)^{\alpha}\|f\|_{B,Q}&=&
l(Q)^{\alpha}\|\chi_{\cup_{i=1}^{N}J_i}f\|_{B,Q}\\
&\leq& l(Q)^{\alpha} \sum_{i=1}^{N}\|\chi_{J_i}f\|_{B,Q}\leq N
\frac{t}{2^d}\le t,
\end{eqnarray*}
contradicting (\ref{para un cubo cualquiera}). Using that $l(Q)\leq
l(J_1)<2l(Q)$ we can also show that

$$
\frac{t}{2^d}<l(Q)^\alpha\|\chi_{J_1}f\|_{B,Q}\leq 2^n\,
l(J_1)^\alpha \|f\|_{B,J_1}
$$
and $Q\subset 3J_1$.
\end{dem}

\medskip

\section{Proof of the main results}

\begin{prueba}{Theorem (\ref{tipo debil de la malphaB})}

Fix a non-negative function $f \in L_\mu^B(\R^d)$. Fix $t>0$ and
define
$$
E_t=\{x\in \R^d:\mathcal{M}_{\alpha,B}f(x)>t\}.
$$
If $t$ is such that the set $E_t$ is empty, we have nothing to prove.
Otherwise, for each $x \in E_t$ there exists a cube $Q_x$ containing $x$ such that
$$
l(Q_x)^{\alpha}\|f\|_{B,Q_x}>t.
$$
By Lemma \ref{paso de un cubo a un diadico}, there exists a constant $\beta$ and a dyadic cube $P_x$ with $Q_x\subset 3 P_x$ such that
\begin{equation}\label{la 33}
l(P_x)^{\alpha}\|f\|_{B,P_x}>\beta t.
\end{equation}
Since $f \in L^B_\mu(\R^d)$, it is not hard to prove that we can
replace the collection $\{P_x\}$ with a
 maximal disjoint subcollection $\{P_j\}$. Each $P_j$ satisfies (\ref{la 33}) and, by our choice of the $Q_x$'s, $E_t\subset \cup_j 3P_j$.
By Lemmas \ref{el 3.12} and \ref{el 3.13},
$$
\phi_1(\mu(E_t))\leq \sum_j \phi_1(\mu(3P_j)).
$$
On the other hand, inequality (\ref{la 33}) implies that, for each $j$,
$$
\frac{1}{l(P_j)^n}\int_{P_j}B\left(\frac{l(P_j)^{\alpha}|f|}{\beta
t}\right) \, d\mu >1,
$$
and then by the definition and properties of $h_B$,
\begin{eqnarray*}
1&<&\frac{1}{l(P_j)^n}\int_{P_j} B\left(\frac{3^\alpha l(P_j)^\alpha
|f(x)|}{3^\alpha \beta t}\right)d\mu(x)\\
&\leq&
\frac{3^nh_B(3^{-\alpha}\beta^{-1})h_B(l(3P_j)^{\alpha})}{l(3P_j)^n}\int_{P_j}B\left(\frac{|f(x)|}{t}\right)d\mu(x)\\
&\leq&
\frac{C}{\phi_1(l(3P_j)^n)}\int_{P_j}B\left(\frac{|f(x)|}{t}\right)d\mu(x).
\end{eqnarray*}
Hence, since the $P_j$'s are disjoint,
\begin{eqnarray*}
\hspace{-1cm}\sum_j \phi_1(\mu(3P_j))&\leq&  \sum_j
\phi_1(l(3P_j)^n)\\&\leq& C \sum_j \int_{P_j}B\left(\frac{|f(x)|}{t}\right)d\mu(x)\\
&\leq & C  \int_{\R^d}B\left(\frac{|f(x)|}{t}\right)d\mu(x).
\end{eqnarray*}

\end{prueba}

\begin{prueba}{Theorem \ref{teo similar al 5.3 de GC martell}}
% The proof follows similar arguments as those in \cite{garcia cuerva martell}.

Without loss of generality, we can assume that $f$ is a non-negative
bounded function with compact support. This guarantees that
$\mathcal{M}_{\alpha, B}f$ is finite $\mu$-almost everywhere. In
fact, $f\in L_\mu^{p_0}(\R^d)$ where $p_0$ is the exponent of the
hypotheses. From Theorem \ref{acotmaximal} we get that
$\mathcal{M}_{\alpha, B}f\in L_\mu^{q_0}(\R^d)$ and thus
$$\mathcal{M}_{\alpha, B}f(x)<\infty\;\;\;\; \text{a.e.} \;\;x \in \R^d.$$

For each $k\in \Z$ let $\Omega_k=\{x\in \R^d:
2^k<\mathcal{M}_{\alpha,B}f(x)\leq2^{k+1}\}$. Thus $$ \R^d=\cup_{k\in
\Z}\Omega_k.$$
Then, for every $k$ and every $x\in \Omega_k$, by the definition of
$\mathcal{M}_{\alpha, B}f$, there is a cube $Q_x^k$ containing $x$,
such that
$$
l(Q^k_x)^\alpha \|f\|_{B,Q^k_x}>2^k,
$$
and, from Lemma \ref{paso de un cubo a un diadico} there
exist a constant $\beta$ and a dyadic cube $P_x^k$ with $Q_x^k
\subset 3P_x^k$ such that
\begin{equation}\label{6}
l(P^k_x)^\alpha \|f\|_{B,P^k_x}>\beta 2^k.
\end{equation}
From the fact that $B$ is submultiplicative and $\mathop{\rm
supp}(f)$ is a compact set, the inequality above allow us to obtain
that
$$
\frac{l(P^k_x)^n}{B(l(P^k_x)^\alpha)}<\int_{P^k_x}B\left(\frac{|f|}{2^k\beta}\right)d\mu
 \leq C \mu({\rm supp} (f))\leq C.
$$
From the hypotheses on $B$ it is easy to check that

$$
C_1\varphi^{-1}\left(\frac{l(P^k_x)^n}{C}\right)\left(\frac{l(P^k_x)^n}{C}\right)^{\alpha/n}
\leq B^{-1}\left(\frac{l(P^k_x)^n}{C}\right)\leq l(P^k_x)^\alpha,
$$
which allow us to conclude that, for each $k$, $l(P^k_x)$ is bounded
by a constant independent of $x$. Then, there is a subcollection of
maximal cubes (and so disjoint) $\{P^k_j\}_j$ such that every
$Q_x^k$ is contained in $3P^k_j$ for some $j$ and, as a consequence,
$\Omega_k\subset_j 3P^k_j $. Next, decompose $\Omega_k$ into the
sets
$$
E_1^k=3P_1^k\cap \Omega_k, \; E^k_2 =\left(3P^k_2\setminus 3P^k_1\right)\cap \Omega_k, \; ....,\; E^k_j=\left(3P^k_j\setminus \cup_{r=1}^{j-1}3P^k_r \right)\cap \Omega_k,....
$$
Then
$$
\R^d=\bigcup_{k\in \Z}\Omega_k=\bigcup_{j,k}E_j^k
$$
and these sets are pairwise disjoint. Let $K$ be a fixed positive
integer which will go to infinity later, and let $\Lambda_K
=\{(j,k)\in \N \times\Z : |k|\leq K\}$. By using that $E^k_j\subset
\Omega_k$ and that the cubes $P^k_j$ satisfy \eqref{6} we obtain
that

\begin{align*}
\mathcal{I}_k &= \int_{\cup_{-K}^{K}\Omega_k} \left(\mathcal{M}_{\alpha,B}f(x)\right)^qu(x)d\mu(x)\\
&
= \sum_{(j,k)\in \Lambda_k} \int_{E^k_j}\left(\mathcal{M}_{\alpha,B}f(x)\right)^qu(x)d\mu(x)\\
&
\leq \sum_{(j,k)\in \Lambda_k} u(E^k_j)2^{(k+1)q}\\
&
\leq C2^q \sum_{(j,k)\in \Lambda_k}u(E^k_j)\left(l(P^k_j)^\alpha\|f\|_{B,P^k_j}\right)^q\\
& \leq  C 2^q \sum_{(j,k)\in
\Lambda_k}u(3P^k_j)\left(l(P^k_j)^\alpha\|fv^{1/p}\|_{C,P^k_j}\|v^{-1/p}\|_{A,P^k_j}\right)^q,
\end{align*}
 where in the last inequality we have used the generalized H\"{o}lder's inequality and the hypothesis on the functions $A$, $B$ and $C$. Now, by applying the hypothesis on the weights we obtain that
 $$
 \mathcal{I}_k\leq C \sum_{(j,k)\in \Lambda_k}l(3P^k_j)^{nq/p}\|fv^{1/p}\|^q_{C,P^k_j}=C\int_{\mathcal{Y}}T_k(fv^{1/p})^q d\nu,
 $$
where $\mathcal{Y}=\N \times \Z$, $\nu$ is de measure in
$\mathcal{Y}$ given by $\nu(j,k)=l(3P^k_j)^{nq/p}$ and, for every
measurable function $h$, the operator $T_k$ is defined by the
expression
$$
T_kh(j,k)=\| \varphi \|_{C,P^k_j} \chi_{\Lambda_k}(j,k).
$$
Then, if we prove that $T_k:L^p(\R^d,\mu)\to
L^q(\mathcal{Y},\nu)$ is bounded independently of $K$, we shall
obtain that
$$
\mathcal{I}_k\leq C \int_{\mathcal{Y}}T_k(fv^{1/p})^q d\nu\leq
C\left(\int_{\R^d}(fv^{1/p})^pd\mu\right)^{q/p}=C\left(\int_{\R^d}f^pvd\mu\right)^{q/p},
$$
and we shall get the desired inequality by doing $K\to\infty$. But
the proof of the boundedness of $T_k$ follows the same arguments as
in Theorem 5.3 in \cite{garcia cuerva martell}, by using now that
the function $C \in B_p$, so we omit it.

\end{prueba}

%..........................................................

\begin{prueba}{Theorem \ref{teo puntual de Malpha a M}}
Let $g(x)=|f(x)|^{p/s}$, then
$$
|f(x)|=g(x)^{s/p+\alpha/n-1}g(x)^{1-\alpha/n}.
$$
Let $x\in \R^d$ and $Q$ be a fixed cube containing $x$. By the
generalized H\"{o}lder's inequality (\ref{desigualdad de holder
general}) and the fact that
$$
g(x)^{(s/p+\alpha/n-1)n/\alpha}=|f|^p,
$$
we get
\begin{eqnarray*}
l(Q)^\alpha \|f\|_{B,Q}&\leq& C \,l(Q)^\alpha
\|g^{1-\alpha/n}\|_{\phi,Q}\|g^{s/p+\alpha/n-1}\|_{n/\alpha,Q}\\&=&
C \;l(Q)^\alpha
\|g\|_{\psi,Q}^{1-\alpha/n}\left(\frac{1}{l(Q)^n}\int_Q |f(y)|^p
d\mu(y)\right)^{\alpha/n}\\ &\leq& C\;
\left[\mathcal{M}_{\psi}(g)(x)\right]^{1-\alpha/n}\|f\|_{L^p(\mu)}^{{p\alpha}/n}.\\
\end{eqnarray*}
\end{prueba}

\begin{prueba}{Theorem \ref{Meta}}
From Theorem \ref{tipo debil de la malphaB} applied to the case
$\alpha=0$ it is easy to check that
$$
\mu(\{y \in \R^d: \mathcal{M}_B f(y)>2t\})\le C
\int_{|f|>t}B(|f|/t)d\mu(t).
$$
Thus, by changing variables and using inequality above we obtain
that
\begin{eqnarray*}
\int_{\R^d} \mathcal{M}_B f(y)^p d\mu(y)&=&C \int_0^\infty t^p
\mu(\{y \in \R^d: \mathcal{M}_B
f(y)>2t\})\frac{dt}{t}\\
&\leq& C \int_{\R^d}\int_0^{|f(y)|} t^p
B\left(\frac{|f(y)|}{t}\right)\frac{dt}{t}d\mu(y)\\
&=& C \left(\int_{\R^d}|f(y)|^p
d\mu(y)\right)\left(\int_{1}^{\infty}\frac{B(s)}{s^p}\frac{ds}{s}\right).
\end{eqnarray*}
Thus, condition $B_p$ allow us to obtain the desired result.
\end{prueba}

\medskip

\begin{prueba}{Theorem \ref{acotmaximal}}
By Theorem \ref{teo puntual de Malpha a M}, if
$1<p<n/\alpha$, we have
\begin{eqnarray*}
\left(\int_{\R^d}\left(\mathcal{M}_{\alpha,B}(f)\right)^qd\mu\right)^{1/q}&\leq&
C\;
\left(\int_{\R^d}\left(\mathcal{M}_{\psi}(|f|^{p/s})^{1-\alpha/n}\|f\|_{L^p(\mu)}^{p\alpha/n}\right)^qd\mu\right)^{1/q}\\
&=& C\,
\|f\|_{L^p(\mu)}^{p\alpha/n}\left(\int_{\R^d}\mathcal{M}_{\psi}(|f|^{p/s})^{s}d\mu\right)^{1/q}.
\end{eqnarray*}
From Proposition $\ref{Bq}$  we have that the function $\psi\in
B_s$. Thus Theorem $\ref{Meta}$ implies that
$\mathcal{M}_{\psi}:L^s(\mu)\to L^s(\mu)$, and thus
$$
\left(\int_{\R^d}\left(\mathcal{M}_{\alpha,B}(f)\right)^qd\mu\right)^{1/q}\leq
C\, \|f\|_{L^p(\mu)}^{p\alpha/n}
\left(\int_{\R^d}(|f|^{p/s})^sd\mu\right)^{1/q} = C\;
\|f\|_{L^p(\mu)}.
$$
On the other hand, if $p=n/\alpha$ and $Q$ is a cube
such that $x\in Q$ we obtain that
\begin{eqnarray*}
l(Q)^\alpha \|f\|_{\eta,Q}&\le& C l(Q)^\alpha
\|\chi_Q\|_{\phi,Q}\|f\|_{n/\alpha,Q}\\
&\le& C \|f\|_{n/\alpha},
\end{eqnarray*}
and thus
\begin{equation*}
\mathcal{M}_{\alpha,B}(f)(x)\le C \|f\|_{n/\alpha}
\end{equation*}
for a.e. $x$ which leads us with the desired result.
\end{prueba}

%.................................................................

%\section{Acknowledgement}
%
%The first author is supported by Instituto de Matem\'atica Aplicada del Litoral (CONICET-UNL) and Departamento de Matem\'atica
%(FIQ-UNL), Santa Fe, Argentina.
%The second author is  supported by Instituto de Matem\'atica Bah\'ia Blanca (CONICET-UNS) and Departamento de Matem\'aticas (UNS), Bah\'ia Blanca, Argentina.
%%
%Finally, the authors would like to thank the referee whose suggestions and comments
%have been very helpful to improve the presentation of the paper.

 %
 %

%
%
%
%
%

\small
\markright{}

\emph{Gladis Pradolini, Instituto de Matem\'atica Aplicada del Litoral (CONICET-UNL), Departamento de Matem\'atica
(FIQ-UNL), 3000 Santa Fe, Argentina}.
e-mail address: gpradolini@santafe-conicet.gov.ar

\emph{Jorgelina Recchi, Instituto de Matem\'atica Bah\'ia Blanca (CONICET-UNS) Departamento De Matem\'aticas, Universidad Nacional Del Sur, 8000 Bah\'ia Blanca, Argentina.}
e-mail address: drecchi@uns.edu.ar
\end{document}